\documentclass[11pt,a4paper,reqno]{amsart}
\usepackage{mathrsfs,amssymb}

\usepackage{pstricks,pst-node,pst-plot}
\usepackage{epic}
\usepackage{autograph}

\newtheorem{theorem}{Theorem}

\newtheorem{corollary}{Corollary}
\newtheorem{definition}{Definition}
\newtheorem{proposition}{Proposition}
\newtheorem{lemma}{Lemma}

\newcommand{\ol}[1]{\overline{#1}}

\newcommand{\res}[3]{#1|^{#2}_{#3}}
\newcommand{\resmin}[4]{\res{#1}{#2}{#3} \setminus #4}

\begin{document}

\title[Free Left and Right Adequate Semigroups]{Free Left and Right Adequate Semigroups}

\maketitle

\begin{center}

    MARK KAMBITES

    \medskip

    School of Mathematics, \ University of Manchester, \\
    Manchester M13 9PL, \ England.

    \medskip

    \texttt{Mark.Kambites@manchester.ac.uk} \\
\end{center}

\begin{abstract}
Recent research of the author has given an explicit geometric description
of free (two-sided) adequate semigroups and monoids, as sets of labelled directed trees under a
natural combinatorial multiplication. In this paper we show that there are
natural embeddings of each free right adequate and free left adequate
semigroup or monoid into the corresponding free adequate semigroup or monoid.
The corresponding classes of trees are easily described and the resulting
geometric representation of free left adequate and free right adequate
semigroups is even easier to understand than that in the two-sided case.
We use it to establish some basic structural properties of free left and
right adequate semigroups and monoids.
\end{abstract}

\section{Introduction}

\textit{Left adequate semigroups} are an important class of semigroups in which
the right cancellation properties of elements in general are reflected
in the right cancellation properties of the idempotent elements. \textit{Right
adequate semigroups} are defined dually, while semigroups which are both left and
right adequate are termed \textit{adequate}. Introduced 
by Fountain \cite{Fountain79}, these classes of semigroups form a natural
generalisation of \textit{inverse semigroups}, and their study is a
key focus of the \textit{York School} of semigroup theory. Left [right]
adequate
semigroups are most naturally viewed  as algebras of signature $(2,1)$,
with the usual multiplication augmented with a unary operation taking
each element to an idempotent sharing its right [left] cancellation properties.
Within the category of $(2,1)$-algebras the left [right] adequate semigroups form
a quasivariety, from which it follows \cite[Proposition~VI.4.5]{Cohn81}
that there exist free left and right adequate semigroups for every cardinality
of generating set.

When studying any class of algebras, it is very helpful to have
an explicit description of the free objects in the class.
Such a description permits one to understand which identities
do and do not hold in the given class, and potentially to express any
member of the class as a collection of equivalence classes of elements
in a free algebra. In the case of inverse semigroups, a description of
the free objects first discovered by
Scheiblich \cite{Scheiblich72} was developed by Munn \cite{Munn74} into an
elegant geometric representation which has been of immense value in the
subsequent development of the subject. The same approach has subsequently been
used to describe the free objects in a number of related classes of semigroups
\cite{Fountain88,Fountain91,Fountain07} and categories \cite{Cockett06};
however, for reasons discussed in \cite{K_freeadequate}, Munn's approach
is not applicable to adequate semigroups.
In \cite{K_freeadequate}, we gave an explicit geometric representation of
the \textit{free adequate semigroup} on a given set, as a collection of (isomorphism
types of) edge-labelled directed trees under a natural multiplication operation.

The focus of this paper is upon the free objects in the quasivarieties of
left adequate and right adequate semigroups. We show that these
embed into the corresponding free adequate semigroups in a natural way,
as the $(2,1)$-algebras generated by the free generators under the
appropriate operations; the resulting representation of these semigroups
is even easier to understand than that for the free adequate semigroup.
These results combine with \cite{K_freeadequate} to yield
a number of results concerning the structure of free left and right
adequate semigroups.

An alternative approach to free left and right adequate semigroups is 
given by recent work of Branco, Gomes and Gould \cite{Branco09}. Their
construction arose from the fact that free left and right adequate
semigroups are \textit{proper} in the sense introduced in \cite{Branco09}.

In addition to this introduction, this article comprises three
sections. In Section~\ref{sec_preliminaries} we briefly recall the
definitions and elementary properties of adequate semigroups, and
the results of \cite{K_freeadequate} concerning free adequate semigroups.
Section~\ref{sec_leftfas} is devoted to the proof that certain subalgebras
of the free adequate semigroup are in fact the free left adequate and
free right adequate semigroups on the given generating set. Finally,
in Section~\ref{sec_remarks} we collect together some remarks on and
corollaries of our main results.

\section{Preliminaries}\label{sec_preliminaries}

In this section we briefly recall the definitions and basic properties
of left, right and two-sided adequate semigroups (more details of which
can be found in \cite{Fountain79}), and also some of
the main definitions and results from \cite{K_freeadequate} characterising
the free adequate semigroup.

Recall that on any semigroup $S$, an equivalence relation $\mathcal{L}^*$ is
defined by $a \mathcal{L}^* b$ if and only if we have $ax = ay \iff bx = by$
for every $x, y \in S^1$. Dually, an equivalence relation $\mathcal{R}^*$ is
defined by $a \mathcal{R}^* b$ if and only if we have $xa = ya \iff xb = yb$
for every $x, y \in S^1$.
A semigroup is called \textit{left adequate} [\textit{right adequate}] if
every $\mathcal{R}^*$-class [respectively, every $\mathcal{L}^*$-class]
contains an idempotent, and the idempotents commute. A semigroup is
\textit{adequate} if it is both left adequate and right adequate. It is
easily seen that in a
left [right] adequate semigroup, each $\mathcal{R}^*$-class
[$\mathcal{L}^*$-class] must contain a \textit{unique} idempotent.
We denote by $x^+$ [respectively, $x^*$] the unique idempotent in the
$\mathcal{R}^*$-class [respectively, $\mathcal{L}^*$-class] of an element
$x$; this idempotent acts as a left [right] identity for
$x$. The unary operations $x \mapsto x^+$ and $x \mapsto x^*$ are of such critical
importance in the theory of adequate [left adequate, right adequate]
semigroups that it is usual to consider these semigroups as algebras of signature
$(2,1,1)$ [or $(2,1)$ for left adequate and right adequate semigroups]
with these operations. In particular, one restricts attention to
morphisms which preserve the $+$ and/or $*$ operations (and hence
coarsen the 
$\mathcal{R}^*$ and $\mathcal{L}^*$ relations) as well as the
multiplication. Similarly, adequate
[left or right adequate] monoids may be viewed as algebras of signature
$(2,1,1,0)$ [$(2,1,0)$] with the identity a distinguished constant symbol.

The following proposition recalls some basic properties of left and right adequate semigroups; these
are well-known and full proofs can be found in \cite{K_freeadequate}.
\begin{proposition}\label{prop_adequatebasics}
Let $S$ be a left adequate [respectively, right adequate] semigroup and
let $a, b, e, f \in S$ with $e$ and $f$ idempotent. Then
\begin{itemize}
\item[(i)] $e^+ = e$ [$e = e^*$];
\item[(ii)] $(ab)^+ = (ab^+)^+$ [$(ab)^* = (a^*b)^*$];
\item[(iii)] $a^+ a = a$ [$a a^* = a$];
\item[(iv)] $ea^+ = (ea)^+$ [$a^* e = (ae)^*]$;
\item[(v)] $a^+(ab)^+ = (ab)^+$ and [$(ab)^* a^* = (ab^*)$];
\item[(vi)] If $ef = f$ then $(ae)^+(af)^+ = (af)^+$ [$(ea)^* (fa)^* = (fa)^*$].
\end{itemize}
\end{proposition}

Recall that an object $F$ in a concrete category $\mathcal{C}$ is called
\textit{free} on a subset $\Sigma \subseteq F$ if every function from $\Sigma$ to an
object $N$ in $\mathcal{C}$ extends uniquely to a morphism from $F$ to $N$.
The subset $\Sigma$ is called a \textit{free generating set} for $F$, and its
cardinality is the \textit{rank} of $F$.

It is easily seen that classes of left and right adequate semigroups
form a \textit{quasivariety}, and it follows from general results (see, for
example, \cite[Proposition~VI.4.5]{Cohn81}) that free left and right adequate
semigroups and monoids exist. Branco, Gomes and Gould \cite{Branco09} have
recently made the first significant progress in the study of these semigroups. The
main aim of the present paper is to give an explicit geometric representation
of them. We begin with a
proposition, the essence of which is that the distinction between semigroups
and monoids is unimportant. The proof is essentially the same as for the
corresponding result in the (two-sided) adequate case, which can be found
in \cite{K_freeadequate}.

\begin{proposition}\label{prop_monoidsemigroup}
Let $\Sigma$ be an alphabet. The free left adequate [free right
adequate] monoid on $\Sigma$ is
isomorphic to the free left adequate [free right adequate] semigroup
on $\Sigma$ with a single adjoined element which is an identity for
multiplication and a fixed point for $+$ [respectively, $*$].
\end{proposition}

We now recall some definitions and key results from \cite{K_freeadequate};
a more detailed exposition may be found in that paper. We are concerned with
\textit{labelled directed trees}, by which we mean edge-labelled directed graphs
whose underlying undirected graphs are trees. If $e$ is an edge in such a
tree, we denote by $\alpha(e)$, $\omega(e)$ and $\lambda(e)$ the vertex at
which $e$ starts, the vertex at which $e$ ends and the label of $e$
respectively.

Let $\Sigma$ be an alphabet. A \textit{$\Sigma$-tree} (or just
a \textit{tree} if the alphabet $\Sigma$ is clear) is a directed
tree with edges labelled by elements of $\Sigma$, and with two distinguished
vertices (the \textit{start} vertex and the \textit{end} vertex) such that
there is a (possibly empty) directed path from the start vertex to the end vertex.
Figure~1 shows some examples of $\Sigma$-trees where $\Sigma = \lbrace a, b \rbrace$;
in each tree, the start and end vertices are marked by an arrow-head and a
cross respectively.

A tree with only one vertex is called \textit{trivial}, while a tree with
start vertex equal to its end vertex is called \textit{idempotent}. A
tree with a single edge and distinct start and end vertices is called a
\textit{base tree}; we identify each base tree with the label of its
edge.
In any tree, the
(necessarily unique) directed path from the start vertex to the end vertex is called
the \textit{trunk} of the tree; the vertices of the graph which lie on the
trunk (including the start and end vertices) are called \textit{trunk
vertices} and the edges which lie on the trunk are called \textit{trunk edges}.
If $X$ is a tree we write $\theta(X)$ for the set of trunk
edges of $X$.

\begin{figure}\label{fig_first}
\begin{picture}(89,30)
\thicklines
\setloopdiam{10}
\Large
\setvertexdiam{1}

\letvertex A=(15,5)    \drawinitialstate[b](A){}
\letvertex B=(5,15)    \drawstate(B){}
\letvertex C=(5,25)    \drawstate(C){}
\letvertex D=(25,15)    \drawstate(D){$\times$}
\drawedge(A,B){$a$}
\drawedge(B,C){$b$}
\drawedge[r](A,D){$a$}

\letvertex E=(44,5)    \drawinitialstate[b](E){}
\letvertex F=(44,15)    \drawstate(F){$\times$}
\letvertex G=(44,25)    \drawstate(G){}
\drawedge(E,F){$a$}
\drawedge(F,G){$b$}

\letvertex H=(70,5)     \drawinitialstate[b](H){}
\letvertex I=(70,15)    \drawstate(I){$\times$}
\letvertex J=(63,25)    \drawstate(J){}
\letvertex K=(77,25)    \drawstate(K){}
\letvertex L=(84,15)    \drawstate(L){}
\drawedge(H,I){$a$}
\drawedge(I,J){$b$}
\drawedge(I,K){$b$}
\drawedge[r](L,K){$b$}

\end{picture}
\caption{Some examples of $\lbrace a, b \rbrace$-trees.}
\end{figure}

A \textit{subtree} of a tree $X$ is a subgraph
of $X$ containing the start and end vertices, the underlying undirected
graph of which is connected.
A \textit{morphism} $\rho : X \to Y$ of $\Sigma$-trees $X$ and $Y$ is a map taking
edges to edges and vertices to vertices, such that $\rho(\alpha(e)) = \alpha(\rho(e))$,
 $\rho(\omega(e)) = \omega(\rho(e))$ and $\lambda(e) = \lambda(\rho(e))$ for
all edges $e$ in $X$, and which maps the start and end vertex of $X$ to the
start and end vertex of $Y$ respectively.
Morphisms have the expected properties that the
composition of two morphisms (where defined) is again a morphism, while the restriction of a morphism to a
subtree is also a morphism. A morphism maps the trunk edges of its domain
bijectively onto the trunk edges of its image.

An \textit{isomorphism} is a morphism which is bijective on both edges and
vertices. The set
of all isomorphism types of $\Sigma$-trees is denoted $UT^1(\Sigma)$
while the set of isomorphism types of non-trivial $\Sigma$-trees is
denoted $UT(\Sigma)$. The set of isomorphism types of idempotent
trees is denoted $UE^1(\Sigma)$, while the set of isomorphism types of
non-trivial idempotent trees is denoted $UE(\Sigma)$.
Much of the time we shall be formally concerned not with trees themselves
but rather with isomorphism types. However, where no
confusion is likely, we shall for the sake of conciseness ignore the
distinction and implicitly identify trees with their respective isomorphism
types.

A \textit{retraction} of a tree $X$ is an idempotent morphism from
$X$ to $X$; its image is a \textit{retract} of $X$. A tree $X$ is called
\textit{pruned} if it does not admit a
non-identity retraction. The set of all isomorphism types of pruned trees
[respectively, non-trivial pruned trees] is denoted $T^1(\Sigma)$
[respectively, $T(\Sigma)$].
Just as with morphisms, it is readily verified that a composition of
retractions (where defined) is a retraction, and the restriction of a retraction
to a subtree is again a retraction.  A key foundational result from
\cite{K_freeadequate} is the following.
\begin{proposition}\label{prop_pruningconfluent}[Confluence of retracts]
For each tree $X$ there is a unique (up to isomorphism) pruned tree which
is a retract of $X$.
\end{proposition}
The unique pruned retract of $X$ is called
the \textit{pruning} of $X$ and denoted $\ol{X}$.

We now define some \textit{unpruned operations} on (isomorphism types of) trees.
If $X, Y \in UT^1(\Sigma)$ then $X \times Y$ is (the isomorphism type of)
the tree obtained by glueing together $X$ and $Y$, identifying the end vertex
of $X$ with the start vertex of $Y$ and keeping all other vertices and
all edges distinct. If $X \in UT^1(\Sigma)$ then $X^{(+)}$ is (the isomorphism
type of) the tree with the same labelled graph and start vertex of $X$,
but with end vertex the start vertex of $X$.
Dually, $X^{(*)}$ is the isomorphism type of the idempotent tree with the
same underlying graph and end
vertex as $X$, but with start vertex the end vertex of $X$.
It was shown in \cite{K_freeadequate} that the unpruned multiplication
operation $\times$ is a
well-defined associative binary operation on $UT^1(\Sigma)$; the (isomorphism
type of the) trivial tree is an identity element for this operation, and
$UT(\Sigma)$ forms a subsemigroup. The maps
$X \mapsto X^{(+)}$ and $X \mapsto X^{(*)}$ are well-defined
idempotent unary operations on $UT^1(\Sigma)$, and the subsemigroup
generated by their images is idempotent and commutative.

We define corresponding \textit{pruned operations} on $T^1(\Sigma)$ by 
$XY = \ol{X \times Y}$, $X^* = \ol{X^{(*)}}$ and $X^+ = \ol{X^{(+)}}$.
These inherit the properties noted above for unpruned operations, and
have the additional property that the images of the $*$ and $+$ operations
are composed entirely of idempotent elements.
We recall some more key results from \cite{K_freeadequate}
\begin{theorem}\label{thm_morphism}
The pruning map
$$UT^1(\Sigma) \to T^1(\Sigma), \ X \mapsto \ol{X}$$
is a surjective $(2,1,1,0)$-morphism from the set of isomorphism types of
$\Sigma$-trees under unpruned multiplication, unpruned $(*)$ and unpruned $(+)$ with distinguished
identity element to
the set of isomorphism types of pruned trees under pruned multiplication,
$*$ and $+$ with distinguished identity element.
\end{theorem}
\begin{theorem}\label{thm_adequate}
$T^1(\Sigma)$ is a free adequate monoid, freely generated by the set
$\Sigma$ of base trees.
\end{theorem}
\begin{corollary}\label{cor_adequate}
Any subset of $T^1(\Sigma)$ closed under the operations of pruned
multiplication and $+$ [respectively, $*$] forms a left adequate
[respectively, right adequate] semigroup under these operations.
\end{corollary}

If $X$ is a tree and $S$ is a
set of non-trunk edges and vertices of $X$ then $X \setminus S$ denotes
the largest subtree of $X$ (recalling that a subtree must be connected
and contain the start and end vertices, and hence the trunk) which does not
contain any vertices or edges from $S$.
If $s$ is a single edge or vertex we write $X \setminus s$ for $X \setminus \lbrace s \rbrace$.
If $u$ and $v$ are vertices of $X$ such that there is a directed path from
$u$ to $v$ then we shall denote by $\res{X}{u}{v}$ the tree which has the
same underlying labelled directed graph as $X$ but start vertex $u$ and end vertex $v$.
If $X$ has start vertex $a$ and end vertex $b$ then we define
$\res{X}{u}{} = \res{X}{u}{b}$ and $\res{X}{}{v} = \res{X}{a}{v}$ where
applicable.

\section{Free Left Adequate Monoids and Semigroups}\label{sec_leftfas}

In \cite{K_freeadequate} we saw that the monoids $T^1(\Sigma)$ and semigroups
$T(\Sigma)$ are precisely the free objects in the quasivarieties of
adequate monoids and semigroups respectively. In this section, we prove
the main results of the present paper by establishing a corresponding result for
left adequate and right adequate monoids and semigroups. The spirit and outline
of the proof are similar to that of \cite{K_freeadequate}, but the technical
details are in places rather different.

\begin{definition}[Left and right adequate trees]
A $\Sigma$-tree $X$ is called \textit{left adequate} if for each vertex
$v$ of $X$ there is a directed path from the start vertex to $v$, or
equivalently, if every non-trunk edge in $X$ is orientated away from the
trunk. The sets of isomorphism types of left adequate $\Sigma$-trees and left adequate
pruned $\Sigma$-trees are denoted $LUT^1(\Sigma)$ and $LT^1(\Sigma)$
respectively.

Dually, a $\Sigma$-tree $X$ is called \textit{right adequate} if
for each vertex $v$ of $X$ there is a directed path from $v$ to the end
vertex, or equivalently, if every non-trunk edge in $X$ is orientated towards
the trunk. The sets of isomorphism types of right adequate $\Sigma$-trees and right
adequate pruned $\Sigma$-trees are denoted $RUT^1(\Sigma)$ and $RT^1(\Sigma)$
respectively.
\end{definition}

Returning to our examples in Figure~1, the left-hand and
middle tree are left adequate, while the right-hand tree is not, because
of the presence of the rightmost edge which is orientated towards the start
vertex. None of the trees shown are right adequate.

From now on we shall work with left adequate trees and left adequate monoids,
but of course duals for all of our results apply to right adequate
trees and right adequate monoids.

\begin{proposition}\label{prop_leftclosed}
The set $LUT^1(\Sigma)$ of left adequate $\Sigma$-trees contains
the trivial tree and the base trees, and is closed under unpruned
multiplication, unpruned $(+)$, and taking retracts.
\end{proposition}
\begin{proof}
It follows immediately from the definitions that the trivial tree and base
trees are left adequate.

Let $X$ and $Y$ be left adequate trees with start vertices $u$ and 
$v$ respectively. Then $u$ is the start vertex of $X \times Y$, and $X 
\times Y$ has a directed path from $u$ to $v$. Now for any vertex
$w \in X \times Y$, either $w$ is a vertex of $X$ or $w$ is a vertex of $Y$.
In the former case, there is a directed path from $u$ to $w$ in $X$, and hence
in $X \times Y$. In the latter case, there is a directed path from $v$ to
$w$ in $Y$, and hence in $X \times Y$, which composed with the path from $u$
to $v$ yields a directed path from $u$ to $w$. Thus, $X \times Y$ is left adequate.

Consider next the tree $X^{(+)}$. This has the same underlying directed
graph as $X$ and the same start vertex, so it is immediate that it is
left adequate.

Finally, let $\pi : X \to Y$ be a retraction with image $Y$ a subtree of
$X$. Now for any vertex $w$ in $Y$
there is a directed path from the start vertex of $X$ to $w$ in $X$;
since $Y$ is a subtree it is connected and has the same start vertex as
$X$, so this must also be a path in $Y$. Thus, $Y$ is left
adequate.
\end{proof}

\begin{proposition}\label{prop_leftgenerators}
The set $LT^1(\Sigma)$ of pruned left adequate trees is generated as a
$(2,1,0)$-algebra (with operations pruned multiplication and pruned $+$
and a distinguished identity element) by the set $\Sigma$ of base trees.
\end{proposition}
\begin{proof}
The proof is similar to the corresponding one in \cite{K_freeadequate},
so we describe it only in outline.
Let $\langle \Sigma \rangle$ denote the $(2,1,0)$-subalgebra of $LT^1(\Sigma)$
generated by $\Sigma$. We show that every left adequate $\Sigma$-tree
is contained
in $\langle \Sigma \rangle$ by induction on number of edges.
The tree with no edges is the
identity element of $LT^1(\Sigma)$ and so by definition is contained in
$\langle \Sigma \rangle$.
Now suppose for induction that $X \in LT^1(\Sigma)$ has at least one edge, and
that every tree in $LT^1(\Sigma)$ with strictly fewer edges lies in
$\langle \Sigma \rangle$.

If $X$ has a trunk edge then let $v_0$ be the start
vertex of $X$, $e$ be the trunk edge incident with $v_0$, $a = \lambda(e)$
and $v_1 = \omega(e)$. Let
$Y = \resmin{X}{v_0}{v_0}{e}$ and $Z = \resmin{X}{v_1}{}{e}$.
Then $Y$ and $Z$ are pruned trees
with strictly fewer edges than $X$, and so by induction lie in
$\langle \Sigma \rangle$. Now clearly from the definitions we have
$Y \times a \times Z = X$, and since $X$ is pruned using
Theorem~\ref{thm_morphism} we have
$$Y a Z = \ol{Y \times a \times Z} = \ol{X} = X$$
so that $X \in \langle \Sigma \rangle$ as required.

If, on the other hand, $X$ has no trunk edges then let $e$ be any edge
incident with the start vertex $v_0$, and suppose $e$ has label
$a$. Since the tree is left adequate, $e$ must be orientated away from
$v_0$; let $v_1 = \omega(e)$.
We define
$Y = \resmin{X}{v_0}{v_0}{e}$
and
$Z = \resmin{X}{v_1}{v_1}{e}$, and a similar argument to that above
shows that $X = Y (a Z)^+$ where $Y, Z \in \langle \Sigma \rangle$,
so that again $X \in \langle \Sigma \rangle$.
\end{proof}

Now suppose $M$ is a left adequate monoid and $\chi : \Sigma \to M$ is a
function. Our objective is to show that there is a unique $(2,1,0)$-morphism
from $LT^1(\Sigma)$ to $M$ which extends $\chi$.
Following the strategy of \cite{K_freeadequate}, we begin by defining a map $\tau$ from 
the set of idempotent left adequate $\Sigma$-trees to the set $E(M)$ of idempotents
in the monoid $M$.
Let $X$ be an idempotent left adequate $\Sigma$-tree with start vertex
$v$. If $X$ has no edges then we define $\tau(X) = 1$. Otherwise, we
define $\tau(X)$ recursively, in terms of the value of $\tau$ on left
adequate trees
with strictly fewer edges than $X$, as follows.
Let $E^+(X)$ be the set of edges in $X$ which start at the start vertex $v$
and define
\begin{align*}
\tau(X) \ = \ \prod_{e \in E^+(X)}
[\chi(\lambda(e)) \tau(\resmin{X}{\omega(e)}{\omega(e)}{e})]^+.
\end{align*}
It is easily seen that each $\resmin{X}{\omega(e)}{\omega(e)}{e}$ is a left
adequate tree with strictly fewer edges than $X$, so this gives a
valid recursive definition of $\tau$. Moreover, the product is non-empty
and because idempotents commute in the left adequate monoid $M$, its value
is idempotent and independent of the order in which the factors are
multiplied. Note that if the left adequate monoid $M$ is in fact adequate
then the function $\tau$ defined here takes the same values
on left adequate trees as the corresponding map defined in \cite{K_freeadequate}.

\begin{proposition}\label{prop_lefttausplit}
Let $X$ be an idempotent left adequate tree with start vertex $v$,
and suppose $X_1$ and $X_2$ are subtrees of $X$ such that $X = X_1 \cup X_2$
and $X_1 \cap X_2 = \lbrace v \rbrace$. Then $\tau(X) = \tau(X_1) \tau(X_2)$.
\end{proposition}
\begin{proof}
Clearly we have $E^+(X) = E^+(X_1) \cup E^+(X_2)$, and for $i \in \lbrace 1, 2 \rbrace$
and $e \in E^+(X_i)$ we have
$$\tau(\resmin{X}{\omega(e)}{\omega(e)}{e}) = \tau(\resmin{X_i}{\omega(e)}{\omega(e)}{e})$$
so it follows that
$$[\chi(\lambda(e)) \tau(\resmin{X}{\omega(e)}{\omega(e)}{e})]^+
= [\chi(\lambda(e)) \tau(\resmin{X_i}{\omega(e)}{\omega(e)}{e})]^+.$$
The claim now follows directly from the definition of $\tau$.
\end{proof}

\begin{corollary}\label{cor_lefttausplit}
Let $X$ be an idempotent left adequate tree with start vertex $v$,
and $e$ an edge incident with $v$. Then
$$\tau(X) = \tau(X \setminus \omega(e)) \ [\chi(\lambda(e)) \tau(\resmin{X}{\omega(e)}{\omega(e)}{e})]^+.$$
\end{corollary}
\begin{proof}
Let $X_1 = X \setminus e = X \setminus \omega(e)$,
let $S$ be the set of edges in $X$ which are incident with $v$
and let $X_2 = X \setminus (S \setminus \lbrace e \rbrace)$ be the maximum
subtree of $X$ containing $e$ but none of the other edges incident
with $v$. Now clearly we have $E^+(X_2) = \lbrace e \rbrace$
so by the definition of $\tau$ we have
$$\tau(X_2) = [\chi(\lambda(e)) \tau(\resmin{X}{\omega(e)}{\omega(e)}{e})]^+.$$
We also have $X = X_1 \cup X_2$ and $X_1 \cap X_2 = \lbrace v \rbrace$
so by Proposition~\ref{prop_lefttausplit}
$$\tau(X) = \tau(X_1) \tau(X_2) = \tau(X \setminus \omega(e))
 [\chi(\lambda(e)) \tau(\resmin{X}{\omega(e)}{\omega(e)}{e})]^+$$
as required.
\end{proof}

Next we define a map $\rho : LUT^1(\Sigma) \to M$, from the set of
isomorphism types of left adequate $\Sigma$-trees to the left adequate monoid $M$.
Suppose a tree $X$ has trunk vertices
$v_0, \dots, v_n$ in sequence. For $1 \leq i \leq n$ let $a_i$ be the label
of the edge from $v_{i-1}$ to $v_i$. For $0 \leq i \leq n$
let $X_i = \resmin{X}{v_i}{v_i}{\theta(X)}$. We define
$$\rho(X) \ = \ \tau(X_0) \ \chi(a_1) \ \tau(X_1) \ \chi(a_2) \ \dots \ \chi(a_{n-1}) \ \tau(X_{n-1}) \ \chi(a_n) \ \tau(X_n).$$
Clearly the value of $\rho$ depends only on the isomorphism type of $X$ so
$\rho$ is indeed a well-defined map from $LUT^1(\Sigma)$ to $M$.
Again, if $M$ is right adequate as well as left adequate then the function
$\rho$ takes the same value on left adequate trees as its counterpart in
\cite{K_freeadequate}.

\begin{proposition}\label{prop_leftrhosplit}
Let $X$ be a left adequate tree with trunk vertices $v_0, \dots, v_n$ in sequence,
where $n \geq 1$. Let $a_1$ be the label of the edge from $v_0$ to $v_1$.
Then
$$\rho(X) = \tau(\resmin{X}{v_0}{v_0}{v_1}) \chi(a_1) \rho(\resmin{X}{v_1}{}{v_0})$$
\end{proposition}
\begin{proof}
Let $X_0, \dots, X_n$ be as in the definition of $\rho$, so that
$$\rho(X) = \tau(X_0) \ \chi(a_1) \ \tau(X_1) \ \chi(a_2) \ \dots \ \chi(a_{n-1}) \ \tau(X_{n-1}) \ \chi(a_n) \ \tau(X_n).$$
It follows straight from the definition that
$$\rho(\resmin{X}{v_1}{}{v_0}) = \tau(X_1) \ \chi(a_2) \ \dots \ \chi(a_{n-1}) \ \tau(X_{n-1}) \ \chi(a_n) \ \tau(X_n)$$
so we have
\begin{align*}
\rho(X) &= \tau(X_0) \ \chi(a_1) \ \rho(\resmin{X}{v_1}{}{v_0}) \\
&= \tau(\resmin{X}{v_0}{v_0}{v_1}) \ \chi(a_1) \ \rho(\resmin{X}{v_1}{}{v_0})
\end{align*}
as required.
\end{proof}

\begin{proposition}\label{prop_leftrhounprunedmorphism}
The map $\rho : LUT^1(\Sigma) \to M$ is a morphism of
$(2,1,0)$-algebras.
\end{proposition}
\begin{proof}
Let $X$ and $Y$ be trees, say with trunk vertices $u_0, \dots, u_m$ and
$v_0, \dots, v_n$ in sequence respectively. For each $1 \leq i \leq m$ let
$a_i$ be the label of the edge from $u_{i-1}$ to $u_i$, and for each
$1 \leq i \leq n$ let $b_i$ be the label of the edge from $v_{i-1}$ to $v_i$.
For each $0 \leq i \leq m$ let
$X_i = \resmin{X}{u_i}{u_i}{\theta(X)}$
and similarly for each $0 \leq i \leq n$ define
$Y_i = \resmin{Y}{v_i}{v_i}{\theta(Y)}$.

Consider now the unpruned product $X \times Y$. It is easily seen that
for $0 \leq i < m$ we have
$$\resmin{(X \times Y)}{u_i}{u_i}{\theta(X \times Y)} = X_i$$
while for $0 < i \leq n$ we have
$$\resmin{(X \times Y)}{v_i}{v_i}{\theta(X \times Y)} = Y_i.$$
Considering now the remaining trunk vertex $u_m = v_0$ of $X \times Y$ we have
$$\resmin{(X \times Y)}{u_m}{u_m}{\theta(X \times Y)}
= \resmin{(X \times Y)}{v_0}{v_0}{\theta(X \times Y)}
= X_m \times Y_0.$$
By Proposition~\ref{prop_lefttausplit} and the definition of unpruned
multiplication we have $\tau(X_m \times Y_0) = \tau(X_m) \tau(Y_0)$.
So using the definition of $\rho$ we have
{\small
\begin{align*}
\rho(X \times Y) &= \tau(X_0) \chi(a_1) \tau(X_1) \dots \chi(a_m) \tau(X_m \times Y_0) \chi(b_1) \tau(Y_1) \chi(b_2) \dots \chi(b_n) \tau(Y_n) \\
&= \tau(X_0) \chi(a_1) \tau(X_1) \dots \chi(a_m) \tau(X_m)  \tau(Y_0) \chi(b_1) \tau(Y_1) \chi(b_2) \dots \chi(b_n) \tau(Y_n) \\
&= \rho(X) \rho(Y).
\end{align*}}

Next we claim that $\rho(X^{(+)}) = \rho(X)^+$. We prove this by induction
on the number of trunk edges in $X$. If $X$ has no trunk edges then
$X = X^{(+)}$ and so using the fact that $\tau(X) \in E(M)$ is fixed by
the $+$ operation in $M$ we have
$$\rho(X^{(+)}) = \rho(X) = \tau(X) = \tau(X)^+ = \rho(X)^+.$$
Now suppose for induction that $X$ has at least one trunk edge and that the claim holds for
trees with strictly fewer trunk edges. Recall that
$$X_0 = \resmin{X}{u_0}{u_0}{\theta(X)} = \resmin{X}{u_0}{u_0}{u_1}$$
and let $Z = \resmin{X}{u_1}{}{u_0}$. Now
\begin{align*}
\rho(X^{(+)}) &= \tau(X^{(+)})  &\text{ (by the definition of $\rho$)} \\
&= \tau(X_0) [\chi(a_1) \tau(Z^{(+)})]^+ &\text{ (by Corollary~\ref{cor_lefttausplit})} \\
&= \tau(X_0) [\chi(a_1) \rho(Z^{(+)})]^+  &\text{ (by the definition of $\rho$)} \\
&= \tau(X_0) [\chi(a_1) \rho(Z)^+]^+  &\text{ (by the inductive hypothesis)} \\
&= \tau(X_0) [\chi(a_1) \rho(Z)]^+ &\text{ (by Proposition~\ref{prop_adequatebasics}(ii))} \\
&= [\tau(X_0) \chi(a_1) \rho(Z)]^+ &\text{ (by Proposition~\ref{prop_adequatebasics}(iv))} \\
&= \rho(X)^+ &\text{ (by Proposition~\ref{prop_leftrhosplit})}
\end{align*}
as required.

Finally, it follows
directly from the definition that $\rho$ maps the identity element in
$LUT^1(\Sigma)$ (that is, the isomorphism type of the trivial tree) to the
identity of $M$, and so is a $(2,1,0)$-morphism.
\end{proof}

So far, we have closely followed the proof strategy from
\cite{K_freeadequate}, but at this point it becomes necessary to diverge.
This is because the arguments employed in the two-sided case
involve operations on trees which do not preserve left adequacy, and hence
use the $*$ operation in the monoid $M$ even when starting
with left adequate trees.
Instead, the following lemma about left adequate trees (which fails for
general trees) allows us to follow an alternative inductive strategy.

\begin{lemma}\label{lemma_leftrestriction}
Let $\mu : X \to Y$ be a morphism of left adequate trees, let $e$ be an edge
in $X$ and let $v$ be a vertex such that there is a directed path from $\omega(e)$ to $v$. Then
$$\mu(\resmin{X}{\omega(e)}{v}{e}) \subseteq \resmin{Y}{\mu(\omega(e))}{\mu(v)}{\mu(e)}.$$
\end{lemma}
\begin{proof}
Let $X' = \resmin{X}{\omega(e)}{v}{e}$ and $Y' = \resmin{Y}{\mu(\omega(e))}{\mu(v)}{\mu(e)}$.
Notice first that the image $\mu(X')$ is connected and contains
$\mu(\omega(e))$. Since the underlying undirected graph of $Y$ is a
tree, this means that $\mu(X')$ is either contained in $Y'$ as required, or contains
the edge $\mu(e)$; suppose for a contradiction that the latter holds, say
$\mu(e) = \mu(f)$ for some edge $f$ in $X'$. Now since $X$ is left adequate,
there must be a directed path from the start vertex to $\alpha(f)$. But again
$e$ is orientated away from start vertex, and $\alpha(f)$ is in $X'$, which
is a connected component of $X$ including $\omega(e)$ but not $e$, so this
path must clearly pass through the edge $e$. Let $P$ denote the suffix of this
path which leads from $\omega(e)$ to $\alpha(f)$. Then $\mu(e P)$ is a non-empty directed path in $Y$ from
$\mu(\alpha(e))$ to $\mu(\alpha(f)) = \mu(\alpha(e))$, which 
contradicts the fact that $Y$ is a directed tree.
\end{proof}

\begin{lemma}\label{lemma_lefttauabsorb}
Suppose $\mu : X \to Y$ is a morphism of idempotent left
adequate trees.
Then $\tau(Y) \tau(X) = \tau(Y)$.
\end{lemma}
\begin{proof}
We use induction on the number of edges in $X$. If $X$ has no
edges then we have $\tau(X) = 1$ so the result is clear.
Now suppose $X$ has at least one edge and for induction that the result holds
for trees $X$ with strictly fewer edges. By the definition of $\tau$ we
have
\begin{align*}
\tau(X) \ = \ \prod_{e \in E^+(X)}
[\chi(\lambda(e)) \tau(\resmin{X}{\omega(e)}{\omega(e)}{e})]^+
\end{align*}
while
\begin{align*}
\tau(Y) \ = \ \prod_{e \in E^+(Y)}
[\chi(\lambda(e)) \tau(\resmin{Y}{\omega(\mu(e))}{\omega(\mu(e))}{\mu(e)})]^+.
\end{align*}
Suppose now that $e \in E^+(X)$. Then since $\mu$ is a morphism, the
edge $\mu(e)$ lies in $E^+(Y)$. We claim that the factor corresponding
to $e$ in the above expression for $\tau(X)$ is absorbed into the
corresponding factor for $\mu(e)$ in the above expression for $\tau(Y)$. 

Let $X' = \resmin{X}{\omega(e)}{\omega(e)}{e}$ and $Y' = \resmin{Y}{\omega(\mu(e))}{\omega(\mu(e))}{\mu(e)}$.
By Lemma~\ref{lemma_leftrestriction}, the morphism $\mu$ restricts
to a morphism $\mu' : X' \to Y'$.
Since
$X'$ has strictly fewer edges than $X$, the inductive hypothesis
tells us that $\tau(X') \tau(Y') = \tau(Y')$.
Now by Proposition~\ref{prop_adequatebasics}(vi) we have
$$[\chi(\lambda(e)) \tau(X')]^+ \ [\chi(\lambda(e)) \tau(Y')]^+ \ = \ [\chi(\lambda(e)) \tau(Y')]^+.$$
as required.
\end{proof}

\begin{corollary}\label{cor_lefttausubgraph}
Let $X$ be a subtree of an idempotent left adequate tree $Y$. Then
$\tau(Y) \tau(X) = \tau(Y)$.
\end{corollary}
\begin{proof}
The embedding of $X$ into $Y$ satisfies the conditions of
Lemma~\ref{lemma_lefttauabsorb}.
\end{proof}

\begin{corollary}\label{cor_lefttauretract}
Let $Y$ be a retract of an idempotent left adequate tree $X$.
Then $\tau(X) = \tau(Y)$.
\end{corollary}
\begin{proof}
Let $\pi : X \to X$ be a retract with image $Y$.
Since $\pi$ is a morphism, Lemma~\ref{lemma_lefttauabsorb} tells
us that $\tau(X) \tau(\pi(X)) = \tau(\pi(X)) = \tau(Y)$.
But since $\pi(X)$ is a subgraph of $X$, Corollary~\ref{cor_lefttausubgraph}
yields $\tau(X) \tau(\pi(X)) = \tau(X)$.
\end{proof}

\begin{lemma}\label{lemma_leftrhosplit}
Let $X$ be a left adequate tree with trunk vertices $v_0, \dots, v_n$ in sequence,
where $n \geq 1$. Let $a_1$ be the label of the edge from $v_0$ to $v_1$.
Then
$$\rho(X) = \tau(\res{X}{v_0}{v_0}) \rho(X).$$
\end{lemma}
\begin{proof}
We use induction on the number of trunk edges in $X$.
Let $X' = \res{X}{v_0}{v_0}$. Clearly if $X$ has
no trunk edges then we have $X = X'$ and from the definition of $\rho$ we have
$\rho(X) = \tau(X')$, so the claim reduces to the fact that $\tau(X')$
is idempotent. Now suppose $X$ has at least one trunk edge and that the claim
holds for $X$ with strictly fewer trunk edges. Let $Y = \resmin{X}{v_1}{v_n}{v_0}$,
let $Y' = \res{Y}{v_1}{v_1}$ and let
$X_0 = \resmin{X}{v_0}{v_0}{v_{1}}$.
Let $a_1$ be the label
of the edge from $v_0$ to $v_1$. By Corollary~\ref{cor_lefttausplit}
we have
$$\tau(X') = [\chi(a_0) \tau(Y')]^+ \tau(X_0).$$
Now by Proposition~\ref{prop_leftrhosplit}
 we deduce that $\rho(X) = \tau(X_0) \chi(a_1) \rho(Y)$.
Also, by the inductive hypothesis we have $\rho(Y) = \tau(Y') \rho(Y)$.
Putting these observations together we have
\begin{align*}
\tau(X') \rho(X) \ 
&= \ \left( [\chi(a_1) \tau(Y')]^+ \ \tau(X_0) \right) \  \left( \tau(X_0) \ \chi(a_1) \ \rho(Y) \right) \\
&= [\chi(a_1) \tau(Y')]^+ \ \tau(X_0) \ \chi(a_1) [\tau(Y') \rho(Y)] \\
&= \tau(X_0) \ [\chi(a_1) \tau(Y')]^+ \ [\chi(a_1) \tau(Y')] \ \rho(Y) \\
&= \tau(X_0) \ \chi(a_1) \ \tau(Y') \ \rho(Y) \\
&= \tau(X_0) \ \chi(a_1) \ \rho(Y) \\
&= \rho(X)
\end{align*}
as required.
\end{proof}

\begin{corollary}\label{cor_leftrhobig}
Let $X$ be a left adequate tree with trunk vertices $v_0, \dots, v_n$ in sequence,
where $n \geq 1$. Let $a_1$ be the label of the edge from $v_0$ to $v_1$.
Then
$$\rho(X) = \tau(\res{X}{v_0}{v_0}) \chi(a_1) \rho(\resmin{X}{v_1}{}{v_0}) = \rho(\resmin{X}{}{v_{n-1}}{v_n}) \ \chi(a_n) \ \tau(\resmin{X}{v_n}{v_n}{v_{n-1}}).$$
\end{corollary}
\begin{proof}
We prove the first equality, the rest of the claim being dual. We have
\begin{align*}
\rho(X) &= \tau(\res{X}{v_0}{v_0}) \rho(X) &\text{ (by Lemma~\ref{lemma_leftrhosplit})} \\
&= \tau(\res{X}{v_0}{v_0}) \tau(\resmin{X}{v_0}{v_0}{v_1}) \chi(a_1) \rho(\resmin{X}{v_1}{}{v_0}) &\text{ (by Proposition~\ref{prop_leftrhosplit})} \\
&= \tau(\res{X}{v_0}{v_0}) \chi(a_1) \rho(\resmin{X}{v_1}{}{v_0}) &\text{ (by Corollary~\ref{cor_lefttausubgraph}).}
\end{align*}
\end{proof}

\begin{proposition}\label{prop_leftrhopruning}
Let $X$ be a left adequate tree. Then $\rho(X) = \rho(\ol{X})$.
\end{proposition}
\begin{proof}
Let $\pi : X \to X$ be a retraction with image $\ol{X}$. Suppose $X$ has trunk
vertices $v_0, \dots, v_n$. For $1 \leq i \leq n$ let $a_i$ be the label
of the edge from $v_{i-1}$ to $v_i$. We prove the claim by induction on
the number of trunk edges in $X$. If $X$ has no trunk edges then by the
definition of $\rho$ and Corollary~\ref{cor_lefttauretract} we have 
$$\rho(X) = \tau(X) = \tau(\pi(X)) = \rho(\pi(X)).$$

Next suppose that $X$ has at least one trunk edge, that is, that $n \geq 1$.
Let $Z = \resmin{X}{v_1}{}{v_0}$.
Then by Lemma~\ref{lemma_leftrestriction} we have
$$\pi(Z) = \pi(\resmin{X}{v_1}{}{v_0}) \subseteq \resmin{\pi(X)}{v_1}{}{v_0}= \resmin{\ol{X}}{v_1}{}{v_0}$$
and, since $\pi$ is idempotent with image $\ol{X}$, the converse inclusion
also holds and we have
\begin{equation}\label{eq_b}
\pi(Z) = \resmin{\ol{X}}{v_1}{}{v_0}.
\end{equation}
Moreover, by
Lemma~\ref{lemma_leftrestriction} again, the retraction $\pi$ restricts to a morphism
$\pi' : Z \to Z$. Clearly this morphism must also be a retraction, and $Z$ has
strictly fewer edges than $X$, so by the inductive hypothesis and
Proposition~\ref{prop_pruningconfluent} we have
\begin{equation}\label{eq_c}
\rho(Z) \ = \ \rho(\ol{Z}) \ = \ \rho(\ol{\pi'(Z)}) \ = \ \rho(\pi'(Z)) \ = \ \rho(\pi(Z)).
\end{equation}
It also follows easily from definitions that
\begin{equation}\label{eq_a}
\pi(\res{X}{v_0}{v_0}) = \res{\ol{X}}{v_0}{v_0}
\end{equation}
Now
\begin{align*}
\rho(X) \ &= \ \tau(\res{X}{v_0}{v_0}) \ \chi(a_1) \ \rho(Z) &\text{ (by Corollary~\ref{cor_leftrhobig})}\\
        &= \ \tau(\res{X}{v_0}{v_0}) \ \chi(a_1) \ \rho(\pi(Z))  &\text{ (by \eqref{eq_c})} \\
        &= \ \tau(\pi(\res{X}{v_0}{v_0})) \ \chi(a_1) \ \rho(\pi(Z))  &\text{ (by Corollary~\ref{cor_lefttauretract})} \\
        &= \ \tau(\res{\ol{X}}{v_0}{v_0}) \ \chi(a_1) \ \rho(\resmin{\ol{X}}{v_1}{}{v_0})  &\text{ (by \eqref{eq_b} and \eqref{eq_a})} \\
        &= \ \rho(\ol{X}) &\text{ (by Corollary~\ref{cor_leftrhobig})}.
\end{align*}
\end{proof}

Now let $\hat\rho : LT^1(\Sigma) \to M$ be the restriction of $\rho$ to
the set of (isomorphism types of) pruned left adequate trees.

\begin{corollary}\label{cor_leftrhoprunedmorphism}
The function $\hat\rho$ is a $(2,1,0)$-morphism from $LT^1(\Sigma)$ (with
pruned operations) to the left adequate monoid $M$.
\end{corollary}
\begin{proof}
For any $X, Y \in LT^1(\Sigma)$ by Theorem~\ref{thm_morphism} and
Propositions~\ref{prop_leftrhounprunedmorphism} and \ref{prop_leftrhopruning} we have 
$$\hat\rho(XY) = \rho(XY) = \rho(\ol{X \times Y}) = \rho(X \times Y) = \rho(X) \rho(Y) = \hat\rho(X) \hat\rho(Y)$$
and similarly
$$\hat\rho(X^+) = \rho(\ol{X^{(+)}}) = \rho(X^{(+)}) = \rho(X)^+ = \hat\rho(X)^+.$$
Finally, that $\hat\rho$ maps the identity of $LT^1(\Sigma)$ to the identity
of $M$ is immediate from the definitions.
\end{proof}

We are now ready to prove the main results of this paper, which give a
concrete description of the free left adequate monoid and free right adequate
monoid on a given generating set.
\begin{theorem}\label{thm_onesidedmonoid}
Let $\Sigma$ be a set. Then $LT^1(\Sigma)$ [$RT^1(\Sigma)$] is a
free object in the quasivariety of left [right] adequate monoids, freely
generated by the set $\Sigma$ of base trees.
\end{theorem}
\begin{proof}
We prove the claim in the left adequate case, the right adequate case
being dual.
By Corollary~\ref{cor_adequate}, $LT^1(\Sigma)$ is a left adequate monoid. Now
for any left adequate monoid $M$ and function $\chi : \Sigma \to M$,
define $\hat\rho : LT^1(\Sigma) \to M$ as above.
By Corollary~\ref{cor_leftrhoprunedmorphism}, $\hat\rho$ is a $(2,1,0)$-morphism, and
it is immediate from the definitions that $\hat\rho(a) = \chi(a)$
for every $a \in \Sigma$, so that $\hat\rho$ extends $\chi$.
Finally, by Proposition~\ref{prop_leftgenerators}, $\Sigma$ is a
$(2,1,0)$-algebra generating set for $LT^1(\Sigma)$; it follows that
the morphism $\hat\rho$ is uniquely determined by its restriction to the
set $\Sigma$ of base trees, and hence is the unique morphism with the
claimed properties.
\end{proof}

Combining with Proposition~\ref{prop_monoidsemigroup} we also obtain
immediately a description of the free left adequate and free right adequate semigroups.

\begin{theorem}\label{thm_onesidedsemigroup}
Let $\Sigma$ be a set. Then the $LT(\Sigma)$ [$RT(\Sigma)$] is a
free object in the quasivariety of left [right] adequate semigroups, freely
generated by the set $\Sigma$ of base trees.
\end{theorem}

We also have the following relationship between free adequate, free left
adequate and free right adequate semigroups and monoids.
\begin{theorem}\label{thm_embedding}
Let $\Sigma$ be a set. The free left adequate semigroup [monoid] on $\Sigma$ and free right adequate semigroup [monoid]
on $\Sigma$ embed
into the free adequate semigroup [monoid] on $\Sigma$ as the $(2,1)$-subalgebras
[$(2,1,0)$-subalgebras] generated
by the free generators under the appropriate operations. Their intersection is
the free semigroup [monoid] on $\Sigma$.
\end{theorem}

\section{Remarks and Consequences}\label{sec_remarks}

In this section we collect together some remarks on and consequences of
the results in Section~\ref{sec_leftfas} and their proofs.

In a left adequate tree, the requirement that there be a path from the
start vertex to every other vertex uniquely determines the orientation
on every edge in the tree. Conversely, every edge-labelled \textit{undirected}
tree with given start and end vertex admits an orientation on the
edges which makes it left adequate. It might superficially seem attractive, then,
to identify elements of $LUT^1(\Sigma)$ with \textit{undirected}
edge-labelled trees with distinguished start and end vertices. However, the
reader may easily convince herself that not every retraction of such a tree
defines a retraction of the corresponding directed tree. So in order to
define pruning and multiplication it would be necessary to reinstate the
orientation on the edges, which negates any advantage in dropping the
orientation in the first place.

The construction in Section~\ref{sec_leftfas} of a morphism from $LT^1(\Sigma)$
to a monoid $M$ depends only on the facts that $M$ is associative with
commuting idempotents, and that the $+$ operation is idempotent
with idempotent and commutative image and satisfies the six properties given in the case of
left adequate semigroups by Proposition~\ref{prop_adequatebasics}. So a free
left adequate semigroup is also free in any class of $(2,1,0)$-algebras
which contains it and satisfies these conditions. This includes in
particular the class of \textit{left Ehresmann semigroups}.

As observed in \cite{K_freeadequate}, the classes of monoids we have
studied can be generalised to give corresponding classes of small categories.
A natural extension of our methods can be used
to describe the free left adequate and free right adequate
category generated by a given directed graph.
Just as in the previous remark, the free left adequate category will also be
the free left Ehresmann
category. Left Ehresmann categories are generalisations of the \textit{restriction categories} studied
by Cockett and Lack \cite{Cockett02}, which in the terminology of semigroup
theory are \textit{weakly
left E-ample} categories \cite{GouldAmpleNotes}. The generalisation of our results
to categories thus relates to our main results in the same way that the description of the
free restriction category on a graph given in \cite{Cockett06}
relates to the descriptions of free left ample monoids given by Fountain,
Gomes and Gould \cite{Fountain91,Fountain07}.

To conclude, we note some properties of free left and right adequate
semigroups and monoids, which are obtained by combining Theorem~\ref{thm_embedding}
with results about free adequate semigroups and monoids which were obtained
in \cite{K_freeadequate}.
First of all, since each finitely generated free left adequate or
free right adequate semigroup embeds into a finitely generated adequate
semigroup we have the following.
\begin{theorem}
The word problem for any finitely generated free left or right adequate
semigroup or monoid is decidable.
\end{theorem}
As in the two-sided case, the exact computational complexity of the word
problem remains unclear and deserves further study.

Recall
that an equivalence relation $\mathcal{J}$
is defined on any semigroup by
$a \mathcal{J} b$ if and only if $a$ and $b$ generate the same principal
two-sided ideal. A semigroup is called \textit{$\mathcal{J}$-trivial} if
no two elements generate the same principal two-sided ideal.
\begin{theorem}\label{thm_jtrivial}
Every free left adequate or free right adequate
semigroup or monoid is $\mathcal{J}$-trivial.
\end{theorem}
\begin{proof}
If distinct left [right] adequate $\Sigma$-trees $X$ and $Y$ are $\mathcal{J}$-related in
$LT^1(\Sigma)$ [$RT^1(\Sigma)$] then they are
$\mathcal{J}$-related in the free adequate monoid $T^1(\Sigma)$; but we
saw in \cite{K_freeadequate} that $T^1(\Sigma)$ is $\mathcal{J}$-trivial.
\end{proof}

\begin{theorem}
No free left adequate or free right
adequate semigroup or monoid on a non-empty set is finitely generated
as a semigroup or monoid.
\end{theorem}
\begin{proof}
We saw in \cite{K_freeadequate} that finite subsets of $T^1(\Sigma)$
generate subsemigroups whose trees have a bound on the maximum distance
of any vertex from the trunk. Since $LT^1(\Sigma)$ and $RT^1(\Sigma)$
are subsemigroups containing trees with vertices arbitrarily far from
the trunk, it follows that they cannot even be contained in finitely
generated subsemigroups of $T^1(\Sigma)$, let alone themselves be
finitely generated.
\end{proof}

\section*{Acknowledgements}

This research was supported by an RCUK Academic Fellowship. The author
would like to thank John Fountain and Victoria Gould for their
encouragement and advice, all the authors of \cite{Branco09} for
allowing him access to their unpublished work and work in progress, and
Mark Lawson for alerting him to the existence of work on free restriction
categories \cite{Cockett06}.

\bibliographystyle{plain}

\begin{thebibliography}{10}

\bibitem{Branco09}
M.~J.~J. Branco, G.~M.~S. Gomes, and V.~A.~R. Gould.
\newblock Structure of left adequate and related monoids.
\newblock Preprint available online at {\tt
  www-users.york.ac.uk/{\char`\~}varg1/adequate.pdf}.

\bibitem{Cockett06}
J.~R.~B. Cockett and Xiuzhan Guo.
\newblock Stable meet semilattice fibrations and free restriction categories.
\newblock {\em Theory Appl. Categ.}, 16:No. 15, 307--341 (electronic), 2006.

\bibitem{Cockett02}
J.~R.~B. Cockett and S.~Lack.
\newblock Restriction categories {I}. {C}ategories of partial maps.
\newblock {\em Theoret. Comput. Sci.}, 270(1-2):223--259, 2002.

\bibitem{Cohn81}
P.~M. Cohn.
\newblock {\em Universal algebra}, volume~6 of {\em Mathematics and its
  Applications}.
\newblock D. Reidel Publishing Co., Dordrecht, second edition, 1981.

\bibitem{Fountain79}
J.~B. Fountain.
\newblock Adequate semigroups.
\newblock {\em Proc. Edinburgh Math. Soc. (2)}, 22(2):113--125, 1979.

\bibitem{Fountain88}
J.~B. Fountain.
\newblock Free right {$h$}-adequate semigroups.
\newblock In {\em Semigroups, theory and applications ({O}berwolfach, 1986)},
  volume 1320 of {\em Lecture Notes in Math.}, pages 97--120. Springer, Berlin,
  1988.

\bibitem{Fountain91}
J.~B. Fountain.
\newblock Free right type {A} semigroups.
\newblock {\em Glasgow Math. J.}, 33(2):135--148, 1991.

\bibitem{Fountain07}
J.~B. Fountain, G.~M.~S. Gomes, and V.~A.~R. Gould.
\newblock Free ample monoids.
\newblock {\em Internat. J. Algebra Comput. (to appear)}, 2007.

\bibitem{GouldAmpleNotes}
V.~A.~R. Gould.
\newblock ({W}eakly) left {E}-ample semigroups.
\newblock Notes available online at {\tt
  www-users.york.ac.uk/{\char`\~}varg1/finitela.ps}.

\bibitem{K_freeadequate}
M.~Kambites.
\newblock Free adequate semigroups.
\newblock {\tt arXiv:0902.0297v3 [math.RA]}, 2009.

\bibitem{Munn74}
W.~D. Munn.
\newblock Free inverse semigroups.
\newblock {\em Proc. London Math. Soc. (3)}, 29:385--404, 1974.

\bibitem{Scheiblich72}
H.~E. Scheiblich.
\newblock Free inverse semigroups.
\newblock {\em Semigroup Forum}, 4:351--359, 1972.

\end{thebibliography}

\def\cprime{$'$} \def\cprime{$'$} \def\cprime{$'$}

\end{document}